\documentclass[12pt]{amsart}
\usepackage{amsmath,amsthm,amsfonts,amssymb}
\usepackage{amsmath, amsthm, amssymb, amscd, amsfonts}
\usepackage[all]{xy}
\usepackage{amssymb, amsthm, amsmath}
\usepackage[english]{babel}
\usepackage{amsfonts}
\newcommand{\ra}{\rightarrow}

\newcommand{\ot}{\otimes}
\newcommand{\mtc}{\mathcal}

\newcommand{\lb}{\label}
\newcommand{\Lam}{\Lambda}

\newcommand{\al}{\alpha}
\newcommand{\eps}{\epsilon}
\newcommand{\bn}{\begin}
\newcommand{\en}{\end}
\newcommand{\sub}{\subsection}

\newcommand{\D}{\Delta}

\newcommand{\rh}{\rightharpoonup}

\numberwithin{equation}{section}

\newcommand{\pp}{\perp}
\newcommand{\dw}{\downarrow}
\newcommand{\uw}{\uparrow}

\newcommand{\ch}{\chi}
\newcommand{\mtr}{\mathrm}

\numberwithin{equation}{section}
\title[Semisimple Hopf algebras]
{On normal Hopf subalgebras of semisimple Hopf algebras }
\author{Sebastian  Burciu}\thanks{Research partially supported by  PN-II-RU-PD-2009 CNCSIS grant 14/28.07.2010}
\address{Inst.\ of Math.\ ``Simion Stoilow" of the Romanian Academy
P.O. Box 1-764, RO-014700, Bucharest, Romania }
\email{sebastian.burciu@imar.ro}

\begin{document}
\begin{abstract}A criterion for subcoalgebras to be invariant under the adjoint action is given generalizing Masuoka's criterion for normal Hopf subalgebras. At the level of characters, the image of the induction functor from a normal Hopf subalgebra is isomorphic to the image of the restriction functor.
\end{abstract}
\maketitle
\section*{Introduction}
The character theory of semisimple Hopf algebras has been intensively developed in the last thirty years. Many results analogous to the classical theory of finite groups were obtained. For example, one of the most important results was Zhu's proof that the character ring of a semisimple Hopf algebra is a semisimple ring \cite{Z}. A systematic study of character theory for coalgebras was started in \cite{Lar} and then continued in \cite{NR}, \cite{NR'}, as well as in other papers. An important updated reference on results on character theory and its applications in classification of semisimple Hopf algebras can be found in \cite{Nssr}. In the classification of semisimple Hopf algebras an important role is played by normal Hopf subalgebras. Recently it was proven that a Hopf subalgebra is normal if and only if it is a depth two subalgebra \cite{KB}. One proof of this result uses the character theory for normal Hopf subalgebras developed in \cite{cos}.

This paper contains new results on the character theory of semisimple Hopf algebras.
We study at the character level, the induction functor from a normal Hopf subalgebra of a semisimple Hopf algebra to the whole Hopf algebra. An algebraic characterization of this functor is obtained in this situation. It is also shown that at the level of characters the images of the induction and restriction maps are isomorphic as algebras. As in the group case a character, is in the image of the restriction map if and only if it is invariant under the adjoint coaction. Our approach uses the commuting pair of modules described in \cite{Zind}. An important role in the proof of the above results is played by the module obtained by inducing the trivial module from $H$ to the Drinfeld double $D(H)$.

The paper is organized as follows. In the first section we recall basic results about semisimple Hopf algebras and its character ring.

In Section \ref{commpair} we study co-normal ideals  of a semisimple Hopf algebra $H$. We show that a two sided ideal of a semisimple Hopf algebra $H$ is co-normal if and only if its character as left module is central in $H^*$. This follows by duality from a description of all subcoalgebras of $H$ that are invariant under the adjoint action. We prove that a subcoalgebra of $H$ is invariant under the adjoint action of $H$ if and only if its character as left (or right) $H$-comodule is central in $H$. This generalizes a well known result of Masuoka stating that a Hopf subalgebra of $H$ is normal if and only if its integral is central in $H$ \cite{Masnr}. The proof uses the Drinfeld double $D(H)$ and the commuting pair of modules from \cite{Zind}.

In the third section we prove few results about the induction functor from a semisimple Hopf subalgebra $K$ of $H$ to $H$. This functor induces a map at the level of character rings $\mtr{ind}^H_K:C(K) \ra C(H)$. Its image is an ideal inside $C(H)$. We also have an algebra map $\mtr{res}^H_K:C(H) \ra C(K)$ induced by the restriction functor. In these settings, Rieffel's equivalence relation from \cite{Rie} is also described in terms of the center of the Hopf algebra. 

In the last section we consider the situation when $K$ is a normal Hopf subalgebra of $H$. We prove that the images of the restriction and induction maps are isomorphic as algebras.  Another description of this image is given at the end of this section.

We work over an algebraically closed field $k$ of characteristic zero. For a vector space $V$ by $|V|$ is denoted the dimension of $V$. Unless otherwise specified, all the modules are supposed to be left modules. The comultiplication, counit and antipode of a Hopf algebra are  denoted by $\Delta$, $\epsilon$ and $S$, respectively. We use Sweedler's notation $\D(x)=\sum x_1\ot x_2$ for all $x\in H$.  All the other notations for Hopf algebras are those used in \cite{Montg}

\section{Preliminaries}\label{prelim}

\sub{Characters of semisimple Hopf algebras}
Throughout this paper $H$ will be a finite dimensional semisimple Hopf algebra over the algebraically closed field $k$. Then $H$ is also cosemisimple and $S^2=\mtr{Id}$ (see \cite{Lard}).  We use the notation $\Lam_H \in H$ for the idempotent integral of $H$ ( $\eps(\Lam_H)=1$) and $t_H \in H^*$ for the idempotent integral of $H^*$ ($t_H(1)=1)$. Denote by $\mtr{Irr}(H)$ the set of irreducible characters of $H$ and let $C(H)$ be the character ring  of $H$. Then $C(H)$ is a semisimple subalgebra of $H^*$ \cite{Z} and $C(H)=\mtr{Cocom}(H^*)$, the space of cocommutative elements of $H^*$. By duality, the character ring of $H^*$ is a semisimple subalgebra of $H$ and under this identification it follows that $C(H^*)=\mtr{Cocom}(H)$.

If $M$ is an $H$-module with character $\chi$ then $M^*$ is also an $H$-module with character $\chi^*=\chi \circ S$. This induces an involution $``\;^*\;":C(H)\ra C(H)$ on $C(H)$.

If $\mtr{Irr}(H)$ is the set of irreducible $H$-modules then from \cite{Montg} it follows that the regular character of $H$ is given by the formula
\begin{equation}
\label{f1}|H|t_{ _H}=\sum_{\ch \in \mtr{Irr}(H) }\chi(1)\chi.
\end{equation}

The dual formula is
\begin{equation}\label{f2}|H|\Lambda_{ _H}=\sum_{d \in \mtr{Irr}(H^*)}\eps(d)d\end{equation}

One also has $t_{ _H}(\Lam_{ _H})=\frac{1}{|H|}$ \cite{Lard}.
There is an associative nondegenerate bilinear form on $C(H)$ given by \bn{equation}(\ch,\;\mu)=\ch\mu(\Lambda_H). \end{equation} It follows that $(\ch,\;\mu)=m_{ _H}(\ch,\; \mu^*)$ where $m_{ _H}$ is the multiplicity form on $C(H)$. For two modules $M$ and $N$ of $H$ one has $m_{ _H}(\ch_M,\;\ch_N)=\mtr{dim}_k\mtr{Hom}_H(M,\;N)$. The pairs $\{\ch\}_{\ch \in \mtr{Irr}(H)}$ and $\{ \ch^*\}_{\ch \in \mtr{Irr}(H)}$ form dual bases with respect to $(\;,\;)$.


Recall that a Hopf subalgebra $K$ of a semisimple Hopf algebra $H$ is normal if $\sum h_1xSh_2 \in K$ for all $h \in H$ and $x \in K$.

\subsection{Properties of Fourier transform $\mtc{F}$}
Let $H$ be a semisimple Hopf algebra. Then $H^*$ is a left and right $H$-module via $(a \rightharpoonup f)(b)=f(ba)$ and $(f \leftharpoonup a)(b)=f(ab)$. Similarly $H$ is an $H^*$ left and right $H$-module via $f \rightharpoonup a=\sum f(a_2)a_1$ and $a \leftharpoonup f= \sum f(a_1)a_2$.

Let $\mtc{F}:H \ra H^*$ be the Fourier transform of $H$ given by $\mtc{F}(a)=a \rh t_H$ for all $a \in H$. It is well known that $\mtc{F}$ is bijective \cite{Montg}. The inverse of $\mtc{F}$ is given by \bn{equation*}\mtc{F}^{-1}(f)=|H|Sf \rightharpoonup \Lambda_{ _H}.\end{equation*}


\bn{lemma}\lb{descprop}
Let $K$ be a Hopf subalgebra of a semisimple Hopf algebra $H$. Then
\bn{equation}\label{desc}
H^*=\mtc{F}(K) \oplus K^{\perp}.
\end{equation}
\end{lemma}

\bn{proof}
Let $i :K \hookrightarrow H$ be the canonical inclusion of $K$ in $H$.  Since $i^*:H^* \ra K^*$ is a surjective Hopf map one has $i^*(t_{ _H})=t_{ _K}$. On the other hand the map $i^*$ is just restriction to $K$.
It follows that $t_{ _H}\dw_K^H=t_{ _K}$. Suppose that $f\in \mtc{F}(K) \cap K^{\perp}$. Then $f=a \rh t_H$  for some $a \in K$. Then $0=f(x)=t_H(ax)=t_K(ax)$ for all $x \in K$. Nondegeneracy of $t_K$ \cite{LSw} implies $a=0$ Thus $f=0$. A dimension argument finishes the proof.
\end{proof}

It is well known that the map $\mtc{F}$ sends the center $\mtc{Z}(H)$ of $H$ into the character ring  $C(H)$ of $H$ \cite{Montg}.

\bn{lemma} \label{commdr}Let $K$ be a Hopf subalgebra of a semisimple Hopf algebra $H$.
Then the following diagram
\bn{equation}\begin{CD}H @> \mathcal{F}_H>>  H^* \\
@A i AA @V Vi^*V \\
K @> \mathcal{F}_K>> K^*. \end{CD}\end{equation}
is commutative where $i$ is the canonical inclusion.
\end{lemma}

\bn{proof}
As before $i^*(t_{ _H})=t_{ _K}$ where $i^*$ is just restriction map from $H$ to $K$. If $a \in K$ then \bn{equation}i^*(\mtc{F}_H(i(a))=i^*(a \rh t_H)= a \rh t_K=\mtc{F}_K(a)\end{equation}
\end{proof}

\section{The commuting pair of modules}\label{commpair}
Let $H$ be a semisimple Hopf algebra. The Drinfeld double $D(H)$ of $H$ is defined as follows: $D(H) \cong H^{*cop} \otimes H$ as coalgebras. The multiplication on $D(H)$ is given by:
\bn{equation}
(g \ot h)(f \ot l)=\sum g(h_1\rightharpoonup f \leftharpoonup S^{-1}h_3)\ot h_2l.
\end{equation}
Its antipode is given by $S(f \ot h)=S(h)S^{-1}(f)$. If $H$ is semisimple over $k$ then $D(H)$ is also semisimple and cosemisimple \cite{Montg}.

Consider the induced module from $H$ to $D(H)$ given by $A_0=D(H)\ot_Hk$.
Then $A_0 \cong H^*$ where the action is given by $
a.f(x)=f(\sum S^{-1}a_2xa_1)$ and $g.f=gf$ for all $a, \;x  \in H$ and $f,\;g \in A^*$. \cite{KSZ}

Via $\mtc{F}$, the module $A_0$ can also be realized on $H$ as following:
$x.a=\sum x_1aS(x_2)$ and $f.x=x \leftharpoonup S^{-1}f=\sum f(S^{-1}x_1)x_2$, see \cite{CW} for example. Thus:

\bn{prop}
The map $\mtc{F}$ is an isomorphism of $D(H)$-modules.
\end{prop}

It was proven by Zhu \cite{Zind} that
\bn{equation}\mtr{End}_{D(H)}(A_0)=C(H)^{op}\end{equation} and the isomorphism is given by $\ch \mapsto R_{\ch}$ where $R_{\ch}$ is right multiplication by $\ch$ on $H^*$. Therefore the map
\bn{equation}\label{dsch1}
C(H) \ra \mtr{End}_{D(H)}(H^*),\;\;\;\; \ch \mapsto R_{S(\ch)}
\end{equation}
is a ring isomorphism. Also the map
\bn{equation}\label{dsch2}
C(H) \ra \mtr{End}_{D(H)}(H),\;\;\;\; \ch \mapsto (S(\ch) \rightharpoonup  \;\;-)
\end{equation}
is a ring isomorphism.

Let $E_1,E_2,\cdots\;,E_s$ be the set of central primitive idempotents of $C(H)$. From the above facts it follows that the homogeneous $D(H)$- components of $H^*$ are given by $H^*E_i$ for all $1 \leq i \leq s$.

Let $\mtr{Irr}(H^*)$ be the set of irreducible characters of $H^*$. For any $d \in \mtr{Irr}(H^*)$ let $\xi_d \in H^*$ the primitive idempotent corresponding to $d$. It can be checked that $\mtc{F}(d)=\frac{1}{\eps(d)}\xi_{d^*}$ where $d^*=S(d)$ (see also \cite{Montg} for the dual version.) It is also known that $\mtc{F}$ sends the algebra $C(H^*)$ to the center $\mtc{Z}(H^*)$ of $H^*$.

For a coalgebra $C$ let $\mtr{Irr}(C)$ be the set of irreducible $H^*$-characters contained in $C$. Then $d_{ _C}=\sum_{ _{d \in \mtr{Irr}(C)}}\eps(d)d$ is the character of $C$ as left (or right) $H$-comodule.

\bn{rem}\label{coalg}
Let $C$ be a subcoalgebra of $H$. Then \bn{equation}\mtc{F}(C)=\oplus_{d \in  \mtr{Irr}(C)}H^*\xi_d.\end{equation} Indeed it is enough to verify this equality for a simple subcoalgebra $C$ with character $d$. This follows since $C=d \leftharpoonup H^*$ and $\mtc{F}$ is an isomorphism of $D(H)$-modules.
\end{rem}

\subsection {Co-normal ideals}
A vector subspace $I$ of $H$ is called co-normal if \bn{equation}\sum Sv_3v_1 \ot v_2 \in H\ot I\end{equation} for all $v \in I$. (usually $I$ will be an ideal.)

Let $I$ be an ideal of $H$ and $\pi :H \ra H/I$ the canonical projection. Then $(H/I)^*$ is a subcoalgebra of $H^*$ via $\pi^*$.

\bn{prop}\label{nrchr}
Let $I$ be an ideal of $H$. Then $I$ is a co-normal ideal if and only if the subcoalgebra $(H/I)^*$ is invariant under the adjoint action $H^*$ on itself.
\end{prop}

\bn{proof}
Note that $(H/I)^*=I^{\pp}$ inside $H^*$. $(H/I)^*$ is invariant under the adjoint action of $H^*$ if and only if $\sum f_1gS(f_2) \in I^{\perp}$ for all $f \in H^*$ and $g \in I^{\perp}$. This is equivalent to $(\sum f_1gS(f_2))(x)=0$ for all $f \in H^*$, $g \in I^{\perp}$ and $x \in I$. But $(\sum f_1gS(f_2))(x)= f( \sum x_1S(x_3)g(x_2))$ which implies that $\sum x_1S(x_3)g(x_2)=0$ for all $g \in I^{\pp}$ and $x \in I$. This is equivalent to $\sum x_1S(x_3) \ot x_2 \in H \ot I$.
\end{proof}

Let $A$ be any finite dimensional semisimple algebra over $k$ and $M$ be an $A$-module. Recall that an $A$-submodule of $M$ is called full isotopic submodule of $M$ if and only if it is a sum of homogeneous components of $M$. It is easy to see that $N$ is a full isotopic submodule of $M$ if and only if it is fixed by any $A$-endomorphism of $M$.

\bn{rem}
Let $A$ be a finite dimensional semisimple algebra over an algebraically closed field of characteristic zero. Let $x \in A$ be an idempotent of $A$ and $e \in \mtc{Z}(A)$ be a central idempotent of $A$. Then

1) $Ax$ is a two sided ideal of $A$ if and only if $x$ is central.

2) $Ax =Ae$ if and only if $x=e$.
\end{rem}

\bn{thm}\label{centity}
A subcoalgebra $C$ of $H$  is invariant under the adjoint action if and only if its character $d_{ _C}$ is central in $H$.
\end{thm}

\bn{proof}

Suppose that $C$ is a subcoalgebra of $H$ invariant under the adjoint action. Then $C$ is a $D(H)$-submodule of $H$. Description \ref{dsch2} of the endomorphism ring of $H$ shows that $C$ is invariant under any $D(H)$-endomorphism of $H$. Thus $C$ is a full isotopic submodule of $H$. Since $\mtc{F}$ is an isomorphism of $D(H)$-modules it follows that $\mtc{F}(C)$ is also a full isotopic submodule of $H^*$. Therefore $\mtc{F}(C)=\oplus_{i \in X} H^*E_i=H^*(\sum_{i \in X}E_i)$ for some set $X$. On the other hand by Remark \ref{coalg} one has $\mtc{F}(C)=H^*(\sum_{ _{d \in \mtr{Irr}(C)}}\xi_d)$. The previous remark implies that $\sum_{ _{d \in \mtr{Irr}(C)}}\xi_d=\sum_{i \in X}E_i$. It follows that $\sum_{ _{d \in \mtr{Irr}(C)}}\xi_d \in C(H)$ and therefore $\sum_{ _{d \in \mtr{Irr}(C)}}\eps(d)d^*=\mtc{F}^{-1}(\sum_{ _{d \in \mtr{Irr}(C)}}\xi_d )$ is central in $H$. But then $d_{ _C}=S(\sum_{ _{d \in \mtr{Irr}(C)}}\eps(d)d^*)$ is also central in $H$.

Conversely suppose that the character $d_{ _C}=\sum_{ _{d \in \mtr{Irr}(C)}}\eps(d)d$ is central in $H$. It follows that $\mtc{F}(d_{ _C}) \in C(H)$. But $\mtc{F}(d_{ _C}) =\sum_{ _{d \in \mtr{Irr}(C)}}\xi_{d^*}$ is a central element in $H^*$ and therefore a central character in $C(H)$. This implies that $\sum_{ _{d \in \mtr{Irr}(C)}}\xi_{d^*}=\sum_{i \in X} E_i$. Remark \ref{coalg} implies that $\mtc{F}(C)=\sum_{i \in X} H^*E_i$. Thus $\mtc{F}(C)$ is a $D(H)$-submodule of $H^*$. Since $\mtc{F}$ is an isomorphism of $D(H)$-modules this implies that also $C$ is a $D(H)$-submodule of $H$. Thus $C$ is invariant under the adjoint action.
\end{proof}


\bn{cor}\label{normideals}
A two sided ideal $I$ of $H$ is co-normal if and only if its character as left (or right) $H$ -module is central in $H^*$.
\end{cor}

\bn{proof}
By Proposition \ref{nrchr} $(H/I)^*$ is a normal subcoalgebra of $H^*$. Theorem \ref{centity} implies the conclusion by duality.
\end{proof}

\section{General results on induction and restriction}\label{genres}
Let $K$ be a Hopf subalgebra of a semisimple Hopf algebra $H$. Then the restriction functor from $H$ to $K$ defines an algebra map at the level of characters ring \bn{equation}\mtr{res}_K^H: C(H) \ra C(K).\end{equation} Similarly, the induction  functor from $K$ to $H$ induces a linear map \bn{equation}\mtr{ind}_K^H:C(K) \ra C(H).\end{equation}

For an irreducible character $\al \in \mtr{Irr}(K)$ let $f_{\al}$ be the central idempotent corresponding to $\al$. Similarly if $\ch \in \mtr{Irr}(H)$ then $e_{\ch}$ is the corresponding central idempotent in $H$.
Consider the commutative algebra $\mtc{Z}(H)\cap K $ as a subalgebra of $\mtc{Z}(H)$ and $\mtc{Z}(K)$. Then there are partition of characters $\mtr{Irr}(H)=\bigsqcup_{i=1}^s\mtc{A}_i$ and $\mtr{Irr}(K)=\bigsqcup_{i=1}^s\mtc{B}_i$ such that a basis of primitive idempotents for the above algebra is given by

\bn{equation}m_i=\sum_{ \ch \in \mtc{A}_i}e_{\ch}=\sum_{ \al \in \mtc{B}_i}f_{\al}.\end{equation}

\bn{prop}\label{div1}
With the above notations

\bn{equation}(\sum_{ \ch \in \mtc{A}_i}\ch(1)\ch)\dw_K^H=\frac{|H|}{|K|}\sum_{ \al \in \mtc{A}_i}\al(1)\al.\end{equation}

\bn{equation}(\sum_{ \al \in \mtc{A}_i}\al(1)\al)\uw^H_K=\sum_{ \ch \in \mtc{A}_i}\ch(1)\ch.\end{equation}
\end{prop}

\bn{proof}
By formula \ref{f1} it follows that \bn{equation}t_{ _H}=\frac{1}{|H|}\sum_{\ch \in \mtr{Irr}(H)}\chi(1)\chi.\end{equation} This implies that $\mtc{F}_H(e_{ _{\ch}})=\frac{\ch(1)}{|H|}\ch$. Similarly $\mtc{F}_K(f_{ _{\al}})=\frac{\al(1)}{|K|}\al$.
The first statement follows from the commutativity of the diagram from Proposition \ref{commdr}. Indeed,

\bn{equation}\frac{1}{|K|}\sum_{ \al \in \mtc{A}_i}\al(1)\al=\mtc{F}_K(m_i)=(\mtc{F}_H(m_i))\dw^H_K=(\frac{1}{|H|}\sum_{ \ch \in \mtc{A}_i}\ch(1)\ch)\dw_K^H.\end{equation}

For the second equality note that \bn{equation}(\sum_{ \al \in \mtr{Irr}(K)}\al(1)\al)\uw^H_K=\sum_{\ch \in \mtr{Irr}(H)}\chi(1)\chi.\end{equation} since both terms are the regular characters of $H$. Frobenius reciprocity applied to the first equality of this Proposition implies the second equality.
\end{proof}

\subsection{Image of the induction map}

The following result is Proposition 2 of \cite{Bd}. It shows that the image of the induction map is a two sided ideal in $C(H)$. For a different proof of that see \cite{BKK}.

\bn{lemma}\label{pr}
Let $K$ be a Hopf algebra of a semisimple Hopf algebra. Let $M$ be an $H$ module and $V$ a $K$-module. Then
\bn{equation}M\ot V\uw^H_K=(M\dw^H_K\ot V)\uw^H_K\end{equation} and
\bn{equation}V\uw^H_K\ot M =(V\ot M\dw^H_K)\uw^H_K\end{equation}
\end{lemma}

Let $\eps_K$ the character of the trivial $K$-module. Let $\eps_{ _K}\uw^H_K$ be the character corresponding to the trivial $K$-module induced to $H$.

\bn{prop}\label{incls}
Let $K$ be a Hopf subalgebra of a  semisimple Hopf algebra $H$. Then:

1)$\eps_{ _K}\uw^H_KC(H) \subset \mtr{Im}(\mtr{ind^H_K})$

2)$\mtc{F}(K)\cap C(H) \subset \mtr{Im}(\mtr{ind^H_K})$.
\end{prop}

\bn{proof}
Put $V=k$, the trivial $K$-module, in the second formula of the above lemma. In terms of the characters this becomes $\ch\eps\uw^H_K=\ch\dw_K^H\uw_K^H$ for all $\ch \in C(H)$. This implies that $\eps\uw^H_KC(H) \subset \mtr{Im}(\mtr{ind^H_K})$.

Using the above notations, since $\mtc{F}(K)\cap C(H)=\mtc{F}(K\cap \mtc{Z}(H))$ and $\mtc{F}$ is bijective it follows that $\mtc{F}(m_i)$ form a basis on $C(H)\cap \mtc{F}(K)$. Proposition \ref{div1} shows that $\mtc{F}(m_i)$ are induced modules.
\end{proof}

\section{Induction functor from normal Hopf subalgebras}\label{normal}

\subsection{Restriction to normal Hopf subalgebras}
Let $H$ be a semisimple Hopf algebra over the algebraically closed field $k$ and let  $K$ be a normal Hopf subalgebra of $H$. Define an equivalence relation on the set  $\mtr{Irr}(H)$  by $\ch \sim \mu$ if and only $m_{ _K}(\ch \dw_K^H,\;\mu\dw_K^H)>0$ (see also \cite{Rie}). This is the equivalence relation $r^{H^*}_{ _{L^*,\;k}}$ from \cite{cos} where $L=H//K$. It is proven that $\ch \sim \mu$ if and only if $\frac{\ch\dw_K^H}{\ch(1)}=\frac{\mu\dw_K^H}{\mu(1)}$ (see Theorem 4.3 of \cite{cos}). Thus the restriction of $\ch$ and $\mu$ to $K$ either have the same irreducible constituents or they do not have common constituents at all. This equivalence relation was also considered in \cite{BKK}  in order to define the depth of a Hopf subalgebra.

The above equivalence relation determines an equivalence relation on the set of irreducible characters of $K$. Two irreducible $K$-characters $\al $ and $\beta$ are equivalent if and only if they are constituents of  $\ch\dw_{ _K}^H$ for some irreducible character $\ch$ of $H$. Let $\mtc{C}_1,\cdots \mtc{C}_{s'}$ be the equivalence classes of the equivalence relation defined on $\mtr{Irr}(H)$. Let $\mtc{D}_1,\cdots \mtc{D}_{s'}$ be the corresponding equivalence classes of the equivalence relation defined on $\mtr{Irr}(K)$.

For each $1 \leq  i \leq s'$ put $\al_i=\sum_{\al\in \mtc{D}_i}\al(1)\al \in C(K)$ and $a_i=\sum_{\ch \in \mtc{C}_i}\ch(1)\ch \in C(H)$.

The formulae from subsection $4.1$ of \cite{cos} can be written as:

\bn{equation} \label{restrform}
\ch\dw_{ _K}^H=\frac{\ch(1)}{\al_i(1)}\al_i
\end{equation}

for all $\ch \in \mtc{C}_i$

and

\bn{equation}\label{indform}
\al\uw_{ _K}^H=\frac{\al(1)}{a_i(1)}\frac{|H|}{|K|}a_i
\end{equation}
for all $\al \in \mtc{D}_i$.

\bn{rem}\lb{comb}
Combining the above two formulae it can easily be seen that $\al\uw_{ _K}^H\dw_K^H\uw^H_K=\frac{|H|}{|K|}\al\uw_{ _K}^H$ for all $\al \in \mtr{Irr}(K)$.
\end{rem}

Next proposition relates the partitions from the previous section and those coming from the above equivalence relation.

\bn{prop}
Let $K$ be a normal Hopf subalgebra of a semisimple Hopf algebra $H$. With the above notations one has

1) $s=s'$,

2) $\{\mtc{C}_1,\cdots, \mtc{C}_s\}=\{\mtc{A}_1,\cdots, \mtc{A}_s\}$,

3) $\{\mtc{D}_1,\cdots, \mtc{D}_{s}\}=\{\mtc{B}_1,\cdots, \mtc{B}_s\}$.
\end{prop}

\bn{proof}
Let \bn{equation}p_i=\sum_{\al \in \mtc{D}_i}f_{ _{\al}}\end{equation} for all $1 \leq i \leq s'$.
First we will show that $p_i=\sum_{ \ch \in \mtc{C}_i}e_{ _{\ch}}$ for $1 \leq i \leq s'$.

Since both terms of the previous equality are idempotents it is enough to verify that any irreducible character of $H$ takes the same value on both terms.
Formula \ref{restrform} shows that $\ch(\sum_{\al \in \mtc{D}_i}f_{ _{\al}})=0$ if $\ch \notin \mtc{C}_i$. On the other hand, if $\ch \in \mtc{C}_i$ then \bn{equation}\ch(\sum_{\al \in \mtc{D}_i}f_{ _{\al}})=\sum_{\al \in \mtc{D}_i}m(\al ,\;\ch\dw_K^H)\al(f_{ _{\al}})=\sum_{\al \in \mtc{D}_i}m(\al ,\;\ch\dw_K^H)\al(1)=\ch(1)\end{equation} and the equality is proved.

Thus each element $p_i=\sum_{\al \in \mtc{D}_i}f_{ _{\al}}$ belongs to $\mtc{Z}(H)\cap K $. Since these elements are idempotent each of them is a sum of some of the primitive idempotents $m_j$ of $\mtc{Z}(H)\cap K $. But Proposition \ref{div1} shows that if $\al \in \mtc{B}_i$ and $\beta \in \mtc{B}_j$ with $i \neq j$ then $\al \nsim \beta$. Thus each element $p_i$ coincide with one of the primitive idempotents
 $m_j$ defined above. Since $\sum_{i=1}^{s'}p_i= \sum_{i=1}^sm_i$ it also follows that any $m_j$ coincide with one of the idempotents $p_i$. This fact implies all three equalities.
\end{proof}

Next Proposition is an improvement of Proposition \ref{incls} in the case of a normal Hopf subalgebra.

\bn{prop}\label{fstdescr}
If $K$ is a normal Hopf subalgebra of a semisimple Hopf algebra $H$ then
$\eps\uw^H_KC(H) = \mtr{Im}(\mtr{ind^H_K})=\mtc{F}(K)\cap C(H)$
\end{prop}

\bn{proof}
By Proposition \ref{div1} it follows that $a_i \in \mtc{F}(K)\cap C(H)$ for all $1 \leq i \leq s$.
Since \bn{equation}\al\uw^H_K=\frac{\al(1)}{a_i(1)}a_i\end{equation} for $\al \in \mtc{C}_i$
this shows that that $\mtr{Im}(\mtr{ind^H_K}) \subset \mtc{F}(K)\cap C(H)$. Then the second item of Proposition \ref{incls} implies $\mtr{Im}(\mtr{ind^H_K})=\mtc{F}(K)\cap C(H)$.

On the other hand using Lemma \ref{pr}  and Remark \ref{comb} one has

\bn{equation}(\al\uw^H_K)\eps_{ _K}\uw_K^H=(\al\uw^H_K\dw_K^H)\uw^H_K=\frac{|H|}{|K|}\al\uw^H_K\end{equation}
which shows that $\mtr{Im}(\mtr{ind^H_K}) \subset \eps\uw^H_KC(H) $. The first item of Proposition \ref{incls} implies $\mtr{Im}(\mtr{ind^H_K})=\eps\uw^H_KC(H)$.

\end{proof}

\bn{prop}\lb{isotdesc}
Let $K$ be a normal Hopf subalgebra of $H$.  Then $\mtc{F}(K)$ and $K^{\perp}$ are full isotopic submodules of $H^*$ and the decomposition \bn{equation*}
H^*=\mtc{F}(K) \oplus K^{\perp}.
\end{equation*} of Proposition \ref{descprop} is a decomposition of $D(H)$-modules.
\end{prop}

\bn{proof}
It can be checked directly that $K$ is a submodule of $H$ realized as before. Since $End_{D(H)}(H)=C(H)$ and $K$ is stabilized by any endomorphism of $H$ it follows that $K$ is a full isotopic submodule $H$. It also can be checked that $K^{\perp}$ is stabilized by any endomorphism of $H^*$.
\end{proof}


\subsubsection{General Lemma}
\bn{lemma}
Suppose $A$ is a finite dimensional semisimple algebra and $M$ is a finite dimensional $A$-module. Moreover suppose that $M=M_1\oplus M_2$ where $M_1$ and $M_2$ are full isotopic submodules of $M$. Then the following are true:

\bn{enumerate}
\item \bn{equation*}
\mtr{End}_A(M) \cong \mtr{End}_A(M_1) \oplus  \mtr{End}_A(M_2)
\end{equation*}
as $k$-algebras
\item
\bn{equation*}
 \mtr{End}_A(M_1) =\{F\in \mtr{End}_A(M)\;|\; F(M_2)=0\}
\end{equation*}
\end{enumerate}
\end{lemma}

\bn{prop}
Let $K$ be a normal Hopf subalgebra of $H$. Then:
\bn{enumerate}
\item
\bn{equation*}
\mtr{End}_{D(H)}(K)=\{\ch \in C(H) \;| \;\sum a_1\ch(Sa_2)\in K \; \;\text{ for all}\;\;  a\in H\;\}
\end{equation*}
\item
\bn{equation*}
\mtr{End}_{D(H)}(K^{\perp})=\{\ch \in C(H) \;| \sum  x_1\ch(S(x_2))=0 \; \;\text{ for all}\;\;  x\in K\;\}
\end{equation*}
\end{enumerate}
\end{prop}

\bn{proof}

1)
The above Lemma for the decomposition of Proposition \ref{isotdesc} implies that \bn{equation}\mtr{End}_{D(H)}(K)=\{\ch \in C(H)\;\; \;\; | K^{\perp}R_{S(\ch)}=0\;\}\end{equation}
Since $(fS(\ch))(a)=f(\sum a_1\ch(S(a_2)))$ it follows that $f(\sum a_1\ch(S(a_2)))=0$ for all $f \in K^{\perp}$. Thus $\sum a_1\ch(S(a_2)) \in {K^{\perp}}^{\perp}=K$

2) From the above Lemma it also follows that $\mtr{End}_{D(H)}(K^{\perp})=\{\ch \in C(H)\;\; \;\;| S(\ch) \rightharpoonup K =0 \}$ which is exactly the set mentioned above.
\end{proof}

Define
\bn{equation}\lb{r1}
C^1_H(K)=\{\ch \in C(H) \;| \;\sum a_1\ch(S(a_2))\in K \; \;\text{ for all}\;\;  a\in H\;\}\end{equation}
and
\bn{equation}
\lb{r2}C^2_H(K)=\{\ch \in C(H) \;| \sum x_1\ch(S(x_2))=0 \; \;\text{ for all}\;\;  x\in K\;\}.
\end{equation}
Then the above lemma implies that
\bn{equation}\label{cdesc}
C(H)=C^1_H(K)\oplus C^2_H(K)
\end{equation} as $k$-algebras.

\bn{rem}
Note that $C^2_H(K)=C(H)\cap K^{\perp}=\mtr{ker}(\mtr{res}^H_K)$. Indeed, if $\ch \in C^2_H(K)$ then it follows $\ch(Sx)=0$ for all $x \in K$. The other inclusion is immediate.
\en{rem}

Recall the definition of $\xi_d$ as the central primitive idempotent of $H^*$ corresponding to the irreducible character $d \in \mtr{irr}(H^*)$.

\bn{lemma}
Let $K$ be a normal Hopf subalgebra of a semisimple Hopf algebra $H$. Then
\bn{equation}\lb{formula}\eps_{ _K}\uw_K^H=\frac{|H|}{|K|}\sum_{d \in \mtr{Irr}(K^*)}\xi_d.\end{equation}
\end{lemma}

\bn{proof}
Since $\Lam_K$ is a central element of $H$ it follows that
\bn{equation}\Lam_K=\sum_{\ch \in \mtc{A}}e_{\ch}\end{equation}
for some subset $A\subset \mtr{Irr}(H)$. Applying $\mtc{F}$ to this equality one gets that:
\bn{equation}\mtc{F}(\Lam_K)=\frac{1}{|H|}\sum_{\ch \in \mtc{A}}\ch(1)\ch.\end{equation}

On the other hand Proposition \ref{div1} implies that
\bn{equation}(\sum_{\ch \in \mtc{A}}\ch(1)\ch)\dw_K^H=\frac{|H|}{|K|}\eps_K.\end{equation}

The same Proposition implies that all the other characters of $H$ that are not contained in $\mtc{A}$ do not contain $\eps_K$ when restricted to $K$. Frobenius reciprocity then gives

\bn{equation}\eps_{ _K}\uw_K^H=\sum_{\ch \in \mtr{Irr}(H)}m_K(\ch\dw_K^H,\;\eps_K)=\sum_{\ch \in \mtc{A}}\ch(1)\ch.\end{equation}

Thus $\eps_{ _K}\uw_K^H=|H|\mtc{F}(\Lam_K)$.

On the other hand, since $\Lam_K=\frac{1}{|K|}\sum_{d \in \mtr{Irr}(K^*)}\eps(d)d$

and $\mtc{F}(d)=\frac{\xi_{d^*}}{\eps(d)}$ one can write

\bn{equation}\mtc{F}(\Lam_K)=\frac{1}{|K|}\sum_{d \in \mtr{Irr}(K^*)}\xi_{d^*}=\frac{1}{|K|}\sum_{d \in \mtr{Irr}(K^*)}\xi_{d}.\end{equation}

Thus \bn{equation}\eps_{ _K}\uw_K^H=\frac{|H|}{|K|}\sum_{d \in \mtr{Irr}(K^*)}\xi_d.\end{equation}
\end{proof}

Since $H$ is cosemisimple it is the sum of its simple subcoalgebras. These subcoalgebras are in bijective correspondence with the irreducible $H^*$-characters \cite{Lar}. Since $k$ is algebraically closed, if $d \in \mtr{Irr}(H^*)$ then its associated subcoalgebra $C_d$ is a matrix coalgebra. This means that $C_d$ has a $k$- basis $\{x^d_{ij}\}_{1\leq i,j \leq q}$  such that
\bn{equation}\lb{comtr}
\D(x^d_{ij})=\sum_{l=1}^qx^d_{il}\ot x^d_{lj}
\end{equation}
and $\eps(x_{ij})=\delta_{i,j}$ for all $1\leq i,j \leq q$. Moreover $\eps(d)=q$ and the irreducible character $d$ is given by $d=\sum_{i=1}^qx^d_{ii}$. Therefore one can write \bn{equation}H=\oplus_{d \in \mtr{Irr}(H^*)}C_d\end{equation} which shows that $x^d_{ij}$ form a $k$-basis for $H$. Clearly $\mtr{Irr}(K^*) \subset \mtr{Irr}(H^*)$ and \bn{equation}H=\oplus_{d \in \mtr{Irr}(K^*)}C_d.\end{equation}

On the other hand, by its definition it follows that $\xi_d(x_{ij}^{d'})=\delta_{d,\;d'}\delta_{i,\;j}$ for all $d,d'\in \mtr{Irr}(H)$.

\bn{rem}\lb{onb}
Note that the condition $\ch \in C_H^1(K)$ is equivalent to $\ch(x^d_{ij})=0$ for all $d \in \mtr{Irr}(H^*) \setminus \mtr{Irr}(K^*)$. This can be seen from formula \ref{comtr}.
\end{rem}

\bn{thm}\label{main}
Let $K$ be a normal Hopf subalgebra of $H$ and $\mtr{ind}^H_K:C(K)\ra C(H)$ be the character map induced by the induction functor. Then \bn{equation}C_H^1(K)=\mtr{im}(\mtr{ind}_K^H).\end{equation}
\end{thm}

\bn{proof}
Formula \ref{formula} together with \ref{comtr} show that $\eps_{ _K}\uw_K^H\ch \in C^1_H(K)$ for all $\ch \in C(H)$. Indeed for a basis element $a=x^d_{uv}$ the identity from the definition \ref{r1} of $C_H^1(K)$ is verified since the element $\sum a_1\ch(S(a_2))$ is zero if $d \notin \mtr{Irr}(K^*)$. Obviously $\sum a_1\ch(S(a_2)) \in K$ if  $d \in \mtr{Irr}(K^*)$. Proposition \ref{fstdescr} implies $\mtr{im}(\mtr{ind}_K^H)=\eps_{ _K}\uw_K^HC(H) \subset C_H^1(K)$. On the other hand it will be checked that if $\ch \in C_H^1(K)$ then $\ch\eps_ { _K}\uw_K^H=\frac{|H|}{|K|}\ch$ which shows the other inclusion. The equality $\ch\eps_{ _K}\uw_K^H=\frac{|H|}{|K|}\ch$ is verified by evaluating both terms on the basis elements $x^d_{ij}$ of $H$. If $d \notin \mtr{Irr}(K^*)$ then both evaluations are zero by Remark \ref{onb}. On the other hand if If $d \in \mtr{Irr}(K^*)$ then the evaluations are equal from the definition of $\xi_d$ and the formula for $\eps_{ _K}\uw^H_K$ given in Lemma \ref{formula}.
\end{proof}

\bn{cor}
Let $K$ be a normal Hopf subalgebra of a semisimple Hopf algebra $H$. Then $\mtr{im}(\mtr{ind}_K^H)$ and $\mtr{im}(\mtr{res}_K^H)$ are isomorphic as $k$-algebras.
\end{cor}

\bn{proof}
Since  $C^2_H(K)=\mtr{ker}(\mtr{res}^H_K)$ it follows that the image of the restriction map is isomorphic to $C(H)/C^2_H(K) \cong C^1_H(K)$ as $k$-algebras.
\end{proof}

\subsection{Primitive idempotents for commutative character rings}

In this subsection we will give another description  for the basis of central primitive idempotents of $C(H)$ considered in Section 2.

\bn{prop} Let $V \oplus W= H$ be a decomposition of $H$ as direct sum
of two Drinfeld submodules of $H$. Define $p_{ _V}$ to be the unique linear
functional on $H$ that equals $\eps_{ _H}$ on $V$ and $0$ on $W$. Similarly define $p_{ _W}$. Then $p_{ _V}$ and $p_{ _W}$ are orthogonal idempotents of the character ring of $H$.
\end{prop}
\bn{proof} Suppose $x \in V$. One has
$p_{ _V}^2(x)=p_{ _V}(x_1)p_{ _V}(x_2)=p_{ _V}(x_1)\eps(x_2)=p_{ _V}(x)$ since $\D(x)\in
H\ot V$. On  the other hand if $x \in W$ then $\D(x)\in H\ot W$ and
therefore $p_{ _V}^2(x)=p_{ _V}(x_1)p_{ _V}(x_2)=0=p_{ _V}(x)$. Thus $p_{ _V}^2=p_{ _V}$. It
remains to show that $p_{ _V}$ is a character. For this it is enough to
show that $p_{ _V}(ab)=p_{ _V}(ba)$ for all $a, b \in H$. It is enough to
show that $p_{ _V}(h_1xSh_2)=\eps(h)p_{ _V}(x)$ for all $h, x \in H$. If $x
\in V$ then
$p_{ _V}(h_1xSh_2)=\eps(h_1xSh_2)=\eps(h)\eps(x)=\eps(h)p_{ _V}(x)$. If $x
\in W$ then $p_{ _V}(h_1xSh_2)=0=\eps(x)p_{ _V}(x)$. It is easy to see that $p_{ _V}p_{ _W}=p_{ _W}p_{ _V}=0$.
\end{proof}

\bn{thm}
Let $H=V_1\oplus V_2\cdots V_s$ be the decomposition of $H$
as sum of homogeneous $D(H)$-modules. Then  $\{p_{ _{V_i}}\}$ is the complete set of central primitive orthogonal idempotents of $C(H)$.
\end{thm}

\bn{proof}
Since $\mtr{End}_{D(H)}(H)=C(H)$ it follows that
$s=\mtr{dim}_kZ(C(H))$. Since $\mtc{F}^{-1}$ is a morphism of $D(H)$ -modules one has that $V_i=S(H^*E_i)\rightharpoonup \Lam$ where  $E_1,\;\cdots E_s$ is the complete set of central idempotents of $C(H)$. Then it can be checked that for any $f \in A^*$ one has $<SE_i, S(fE_i)\rightharpoonup \Lam>=\eps_H(S(fE_i)\rightharpoonup \Lam)$ and $<SE_i,S(fE_j)\rightharpoonup \Lam>=0$ for $i \neq j$. This shows that $SE_i=p_{ _{V_i}}$.
\end{proof}

Next we will compute the image of these primitive idempotents under the restriction map $\mtr{res}^H_K : C(H)\ra C(K)$. Suppose that $K=\oplus_{i=1}^{s'}W_i$ is a decomposition of $K$ in homogeneous components as $D(K)$-module. Denote the corresponding central idempotents of $C(K)$ as given in the previous Theorem by $q_{ _{W_i}}$. Then we have the following:

\bn{prop}
Let $K$ be a normal Hopf subalgebra of a semisimple Hopf algebra $H$ and $V\subset H$ a $D(H)$-submodule of $H$.

1) If $V\cap K =0$ then $\mtr{res}_K^H(p_{ _V})=0$.

2) If $V\subset K$ then $V$ can also be regarded as $D(K)$-module and one may suppose that  $V=\oplus_{i=1}^rW_i$ as $D(K)$ module. Then \bn{equation}\mtr{res}^H_K(p_{ _V})=\oplus_{i=1}^rq_{ _{W_i}}.\end{equation}
\end{prop}

\bn{proof}
Note that since $V$ is simple as $D(H)$-module it follows from \ref{isotdesc} that either $V \subset K$ or $V\cap K =0$. The rest of the proof is straightforward.
\end{proof}

{\bf Examples:} 1. Consider $H=kG$. Then it is known that the simple $D(H)$-submodules of $H$ are of the form $k\mathcal{C}:=k<g\;|\;g \in \mathcal{C}>$ where $\mtc{C}$ runs through all conjugacy classes of $G$. Thus in this case $p_{ _V}=p_{ _{\mathcal{C}}}$ is the usual characteristic class function on $\mathcal{C}$.

Let now $K= kN$ with $N$ a normal subgroup of $G$. Clearly if $\mtc{C}\cap N \neq 0$ then $\mtc{C}\subset N$ and $\mtc{C}$ is a union of conjugacy classes of $N$ in this case.

The character ring decomposition from \ref{cdesc} becomes the following in this case:  $$C^1_H(K)=<p_{ _{\mtc{C}}}\;|\;\mtc{C}\subset N>$$ and $$C^2_H(K)=<p_{ _{\mtc{C}}}\;|\;\mtc{C}\cap N= \emptyset >.$$

2. Let $H$ be a semisimple Hopf algebra with $C(H)$ commutative and let $K$ be a normal Hopf subalgebra of $H$ with $C(K)$ also commutative. Using the above Proposition, it follows that the character ring decomposition from \ref{cdesc} is similar to the previous one, namely:

$$C^1_H(K)=<p_{V}\;|\;V\subset K>$$ and $$C^2_H(K)=<p_{V}\;|\;V\cap K=0>.$$
\subsection{Image of the restriction map}

It will be shown that a character of $K$ is in the image of the restriction map if and only if it is invariant under the adjoint coaction. This generalizes a well known result from group algebras.

\bn{thm}\label{descr}
Let $K$ be a normal Hopf subalgebra of a semisimple Hopf algebra $H$.
The image of the restriction map $\mtr{res}^H_K$ is isomorphic as $k$-algebras with the subalgebra of $C(K)$ spanned by \bn{equation}\{ f \in K^*| f(\sum a_1xS(a_2))=\eps(a)f(x)\;\; \text{for all} \;\; a\in A ,x \in K\}.\end{equation}
\end{thm}

\bn{proof}
Clearly the image of the restriction lies in the space \bn{equation*}\{ f \in K^*| f( \sum a_1xS(a_2))=\eps(a)f(x)\;\; \text{for all} \;\; a\in A ,x \in K\}.\end{equation*} Let $v$ be a linear functional of this space. It follows that $v \rightharpoonup $ is an endomorphism of $K$ as $D(H)$-module. Indeed for all $a\in H$ and $x \in K$ one has that:
\bn{equation*}v\rh (a.x)=\sum v \rh a_1xS(a_2)= \sum a_1x_1S(a_4)v(a_2x_2S(a_3)=\end{equation*}\bn{equation*}=\sum a_1x_1S(a_2)v(x_2)=a.(v \rh x)\end{equation*}
Also it can be checked that $v \rh (f.x)=f.(v \rh x)$ for all $x \in K$ and $f \in H^*$. Indeed, \bn{equation*}v \rh (f.x)=v \rh \sum f(Sx_1)x_2=\sum f(Sx_1)x_2v(x_3).\end{equation*} On the other hand, \bn{equation*}f.(v \rh x)=\sum f.v(x_2)x_1=\sum f(Sx_1)x_2v(x_3).\end{equation*}

Since $\mtr{End}_{D(H)}(H)=C(H)$ and $K$ is a a full isotopic $D(H)$-submodule of $H$ it follows that there is $\ch \in C(H)$ such that  $v\rh \;=\ch \rh \;$ on $K$. This implies that $v=\ch\dw_K^H$.
\end{proof}

{\bf Examples:} 1. We start again with the example of normal subgroups. Let $G$ be a finite group and $N$ a normal subgroup of $G$. If $V$ is an
irreducible $N$-module then
\bn{equation*}
V \uw_{ _N}^G\dw_{ _N}^G={\oplus_{i=1}^s}^{g_i}V
\end{equation*}
where $^{g}V$ is a conjugate module of $V$ and $\{g_i\}_{i=1,\;s}$ is a set of representatives
for the left cosets of $N$ in $G$. For $g \in G$ the $H$-module $^gV$
has the same underlying vector space as $V$ and the multiplication with
$h \in N$ is given by $h.v = (g^{-1}hg)v$ for all $v \in V$. 

If $V$ and $V'$ are equivalent modules under the equivalence relation described in the previous section then their restriction to $N$ are scalar multiple one of the other.

Thus the characters of the modules ${\oplus_{i=1}^s}^{g_i}V$ form a basis for the restriction image $\mtr{res}_H^G: C(G) \ra C(N)$ when $V$ runs through the representatives of all the equivalence classes of the above equivalence relation. 

 2. Consider an exact sequence of Hopf algebras:
\bn{equation}\begin{CD}
k @ > >> K  @> i >>  H @ > \pi >> kG@ > >> k 
\end{CD}\end{equation} 
where $G$ is a finite group.

Then $G$ acts on the irreducible characters of $K$ as following. Let $g \in G$ and $h \in H$ with $\pi(h)=g$. If $\ch \in \mtr{Irr}(K)$ then define the character 
\bn{equation}
^g\ch(x):=\ch(h_1xS(h_2))
\end{equation} 
for all $x \in K$. This actions is also described in other papers in the literature. For example, a description at the level of modules is given in \cite{MW}, page 3.

Then the previous Theorem implies that the characters $\sum_{g \in G}{\;\;}^g\ch$ form a basis for the image of the restriction map $\mtr{res}^H_K$ when $\ch$ runs through runs through all the representatives of the classes of the above equivalence relation.

\bibliographystyle{amsplain}
\bibliography{charthv3}

\providecommand{\bysame}{\leavevmode\hbox to3em{\hrulefill}\thinspace}
\providecommand{\MR}{\relax\ifhmode\unskip\space\fi MR }
\providecommand{\MRhref}[2]{%
  \href{http://www.ams.org/mathscinet-getitem?mr=#1}{#2}
}
\providecommand{\href}[2]{#2}
\begin{thebibliography}{10}

\bibitem{Bd}
S.~Burciu, \emph{\textnormal{On some representations of the Drinfel'd double}},
  J. Alg. \textbf{296} (2006), 480–504.

\bibitem{cos}
\bysame, \emph{\textnormal{Coset decomposition for semisimple Hopf algebras}},
  Comm. Alg., \textbf{37} (2009), no.~10, 3573 -- 3585.

\bibitem{KB}
S.~Burciu and L.~Kadison, \emph{\textnormal{Semisimple Hopf algebras and their
  depth two Hopf subalgebras}}, J. Alg. \textbf{1} (2009), 162--176.

\bibitem{BKK}
S.~Burciu, L.~Kadison, and B.~Kuelshammer, \emph{\textnormal{On subgroup
  depth}}, to appear, IEJA (2009), arxiv: 0906.0440.

\bibitem{CW}
M.~Cohen and S.~Westreich, \emph{\textnormal{Some interrelations between Hopf
  algebras and their duals}}, J. Algebra \textbf{283} (2005), 42--62.

\bibitem{KSZ}
Y.~Kashina, Y.~Sommerh\"{a}user, and Y.~Zhu, \emph{\textnormal{Higher
  Frobenius-Schur indicators}}, vol. 181, Mem. Am. Math. Soc., Am. Math. Soc.,
  Providence, RI, 2006.

\bibitem{Lar}
R.~G. Larson, \emph{\textnormal{Characters of Hopf algebras}}, J. Algebra
  (1971), 352--368.

\bibitem{Lard}
R.~G. Larson and D.~E. Radford, \emph{\textnormal{Finite dimensional
  cosemisimple \textnormal{Hopf} Algebras in characteristic zero are
  semisimple}}, J. Algebra \textbf{117} (1988), 267--289.

\bibitem{LSw}
R.~G. Larson and M.~E. Sweedler, \emph{\textnormal{An associative orthogonal
  bilinear form for Hopf algebras}}, Amer. J. Math. \textbf{91} (1969), no.~7,
  75–93.

\bibitem{Masnr}
A.~Masuoka, \emph{\textnormal{Semisimple Hopf algebras of dimension 2p}}, Comm.
  Algebra \textbf{23} (1995), no.~5, 1931--1940.

\bibitem{Montg}
S.~Montgomery, \emph{\textnormal{Hopf algebras and their actions on rings}},
  vol.~82, 2nd revised printing, Reg. Conf. Ser. Math, Am. Math. Soc,
  Providence, 1997.

\bibitem{MW}
S.~Montgomery and S.~Witherspoon, \emph{\textnormal{Irreducible representations
  of crossed products}}, J. Pure Appl. Algebra \textbf{129} (1998), 315--326.

\bibitem{Nssr}
S.~Natale, \emph{\textnormal{Semisolvability of semisimple Hopf algebras of low
  dimension}}, no. 186, Mem. Am. Math. Soc., Am. Math. Soc., Providence, RI,
  2007.

\bibitem{NR}
W.~D. Nichols and M.~B. Richmond, \emph{\textnormal{The Grothendieck group of a
  Hopf algebra}}, J. Pure and Appl. Algebra \textbf{106} (1996), 297--306.

\bibitem{NR'}
\bysame, \emph{\textnormal{The Grothendieck group of a Hopf algebra, I}}, Comm.
  Algebra \textbf{26} (1998), 1081--1095.

\bibitem{Rie}
M.~Rieffel, \emph{\textnormal{Normal subrings and induced representations}}, J.
  Alg. \textbf{24} (1979), 364--386.

\bibitem{Z}
Y.~Zhu, \emph{\textnormal{Hopf algebras of prime dimension.}}, Int. Math. Res.
  Not. \textbf{1} (1994), 53--59.

\bibitem{Zind}
\bysame, \emph{\textnormal{A commuting pair in Hopf algebras}}, Proc. Am. Math.
  Soc. \textbf{125} (1997), 2847--2851.

\end{thebibliography}
\end{document}